\documentclass[preprint,sort&compress,final,12pt]{elsarticle}

\usepackage{amsmath}
\usepackage{amssymb}
\usepackage{graphicx}
\usepackage{latexsym}
\usepackage{epstopdf}
\usepackage{natbib}
\usepackage{color}
\usepackage{hyperref}
\usepackage{amsmath,amscd}

\newtheorem{theorem}{Theorem}[section]

\newtheorem{lemma}[theorem]{Lemma}
\newtheorem{proposition}[theorem]{Proposition}
\newtheorem{remark}[theorem]{Remark}
\numberwithin{equation}{section}

\def\Proof{\noindent{\bf Proof.}~}
\def\qed{\hfill$\square$\smallskip}
\def\dint{\displaystyle\int}

\def\dsup{\displaystyle\sup}

\def\dsum{\displaystyle\sum}
\def\dint{\displaystyle\int}
\def\dsup{\displaystyle\sup}

\def\Im{\mathrm{Im}}

\journal{\empty}
\date{}

\begin{document}

\begin{frontmatter}

\title{The large twist theorem and boundedness of solutions for polynomial potentials with $C^1$ time dependent coefficients}

\author[au1]{Xiong Li\footnote{Corresponding author. Partially supported by the NSFC (11571041) and the Fundamental Research Funds for the Central Universities.}}

\address[au1]{School of Mathematical Sciences, Beijing Normal University, Beijing 100875, P.R. China.}

\ead[au1]{xli@bnu.edu.cn}

\author[au2]{Bin Liu\footnote{Partially supported by the NSFC (11231001).}}

\address[au2]{School of Mathematical Sciences, Peking University, Beijing 100871, P.R. China.}

\ead[au2]{bliu@pku.edu.cn}

\author[au1,au3]{Yanmei Sun}

\address[au3]{School of Mathematics and Information Sciences, Weifang University, Weifang, Shandong, 261061, P.R. China.}

\ead[au3]{sunyanmei2009@126.com}

\begin{abstract}
In this paper we first prove the so-called large twist theorem, then using it to prove the boundedness of all solutions and the existence of quasi-periodic solutions for Duffing's equation
$$
\ddot{x}+x^{2n+1}+\dsum_{i=0}^{2n}p_i(t)x^i=0,
$$
where $p_i(t)\in C^1(\mathbb{S}) (n+1\leq i\leq 2n)$ and  $p_i(t)\in C^0(\mathbb{S}) (0\leq i\leq n)$ with $\mathbb{S}=\mathbb{R}/\mathbb{Z}$.
\end{abstract}


\end{frontmatter}

\section{Introduction}
In this paper we initially are concerned with the boundedness of all solutions and the existence of quasi-periodic solutions for Duffing's equation
\begin{equation}\label{Duffing}
\ddot{x}+x^{2n+1}+\dsum_{i=0}^{2n}p_i(t)x^i=0
\end{equation}
under the smoothness assumption that $p_i(t)\in C^1(\mathbb{S}) \,(n+1\leq i\leq 2n)$ and  $p_i(t)\in C^0(\mathbb{S})\, (0\leq i\leq n)$ with $\mathbb{S}=\mathbb{R}/\mathbb{Z}$.

In the early 1960's, Littlewood \cite{Littlewood} asked whether all solutions of  the general Duffing's
equations
\begin{equation}\label{V}
\ddot{x}+V_x(x,t)=0
\end{equation}
are bounded for all time, that is, $\dsup_{t\in \mathbb{R}} (|x(t)| + \dot{x}(t)|) < +\infty$ holds for all
solutions of Eq.(\ref{V}).

For the Littlewood boundedness problem, during the past years, people have
paid more attention to the special equation (\ref{Duffing}), since
$$
\ddot{x} + x^{2n+1} = 0
$$
is a very nice time-independent integrable system, of which all the solutions are
periodic. Thus if $|x|$ is large enough, Eq.(\ref{Duffing}) can be treated as a perturbation
of an integrable system, then Moser's twist theorem could be applied to prove
the boundedness of all solutions.

The first result was due to Morris \cite{Morris}, who proved that all solutions of a
biquadratic potential
$$
\ddot{x} + 2x^3 = p(t) = p(t + 1)
$$
are bounded. It is noted that here $p(t)$ is only required to be piecewise continuous.

Using the famous Moser's twist theorem \cite{Moser}, Diecherhoff and Zehnder \cite{Diecherhoff}
generalized Morris's result to Eq.(\ref{Duffing}). In that paper, the coefficients $p_i(t)$ are
required to be sufficiently smooth to construct a series change of variables to
transform Eq.(\ref{Duffing}) into a nearly integrable system for large energies. In fact,
in \cite{Diecherhoff}, the smoothness on $p_i(t)$ depends on the index $i$.

An interesting problem (proposed by Diecherhoff and Zehnder) is whether or
not the boundedness of all solutions depends on the smoothness of the coefficients.

In \cite{Laederich}, Laederich and Levi weakened the smoothness requirement on $p_i(t)$ to
$C^{5+\epsilon}$. By modifying the proofs in \cite{Diecherhoff} and using some approximation techniques, the
second author \cite{Liu} obtained the same result for
$$
\ddot{x} + x^{2n+1}+p_1(t)x+p_0(t) = 0,
$$
where the periodic functions $p_1,p_0$ are only required to be continuous, which
shows that the boundedness of all solutions does not depend on the smoothness
of coefficients of lower order terms. Later, Yuan \cite{Yuan95}, \cite{Yuan98}, \cite{Yuan00}
proved that all solutions of Eq.(\ref{Duffing}) are bounded if $p_i(t)\in C^2\,(n+
1 \leq i \leq 2n)$ and $p_i(t) \in C^1 \,(0 \leq i \leq n)$.

Then the remain problem is whether or not the smoothness requirement for
coefficients of higher order terms plays the same role as that for coefficients of
lower order terms.

In \cite{Levi}, Levi and You proved that the equation
$$
\ddot{x} + x^{2n+1} + p(t)x^{2i+1} = 0
$$
with a special discontinuous coefficient $p(t) = K^{[t]mod2}, 0 < K < 1, 2n + 1 >
2i + 1 \geq  n + 2$, possesses an oscillatory unbounded solution. In 2000, Wang \cite{Wang00}
constructed a continuous periodic function $p(t)$ such that the corresponding
equation
$$
\ddot{x} + x^{2n+1} + p(t)x^{i} = 0
$$
possesses a solution which escapes to infinity
in finite time, where $n \geq 2$ and $2n + 1 > i\geq n + 2$. The example constructed in \cite{Wang00} shows that the boundedness of all solutions is
linked with the smoothness of coefficients of higher order terms. Also Wang in \cite{Wang00} pointed out that the continuous periodic function $p(t)$ in this example is not Lipschitz continuous, and whether or not the smoothness requirement in \cite{Wang96}, \cite{Yuan} is the sharpest is unclear.

The aim of this paper is to answer this problem, and prove that all solutions
of Eq.(\ref{Duffing}) are bounded if $p_i(t)\in C^1 \, (n+
1 \leq i \leq 2n)$ and $p_i(t) \in C^0 \, (0 \leq i \leq n)$, which
is the optimal smoothness requirement according to the example constructed in
\cite{Wang00}.

More precisely, we will prove

\begin{theorem}\label{B} If $p_i(t)\in C^1(\mathbb{S}) (n+1\leq i\leq 2n)$ and  $p_i(t)\in C^0(\mathbb{S}) (0\leq i\leq n)$ with $\mathbb{S}=\mathbb{R}/\mathbb{Z}$, then all solutions
of Eq.(\ref{Duffing}) are bounded, that is, every solution of
Eq.(\ref{Duffing}) exists for all $t\in \mathbb{R}$ , and
$$
\sup_{t\in\mathbb{R}}(|x(t)| + |\dot{x}(t)|) < +\infty.
$$
Moreover, there are infinitely many quasi-periodic solutions to Eq.(\ref{Duffing}).
\end{theorem}

In order to prove Theorem \ref{B}, we need to develop Moser's twist theorem and establish the so-called large twist theorem, which will be done in Section 2. And in Section 3, we will prove Theorem \ref{B}.

\section{The large twist theorem}

Consider the mapping
\begin{equation}\label{M}
\mathfrak{M}:\quad \begin{array}{ll}
\left\{\begin{array}{ll}
x_1=x+\beta +\gamma^{-\nu}\,y+ f(x,y;\gamma),\\[0.2cm]
y_1=y+g(x,y;\gamma),\\[0.1cm]
 \end{array}\right.\  \ \ \ \ (x,y)\in \mathbb{R} \times [a,b],
\end{array}
\end{equation}
where $b-a\geq 1$, $\nu\geq 0$ and $\beta$ are two constants, $0<\gamma\leq 1$ is a sufficiently small parameter. Here we assume that $f,g$ are real analytic in $x,y$, continuous in $\gamma$, and have period $2\pi$ in $x$, which can be extended to a complex domain
\begin{equation}\label{D1}
D:\ \  |\Im x|<r_0, \quad y\in D'
\end{equation}
with $D'$ a complex neighborhood of the interval $a\leq y\leq b$, $0<r_0<1$, and the mapping $\mathfrak{M}$ has the intersection property that any curve $\Gamma: y=\phi(x)=\phi(x+2\pi)$ always intersects its image curve $\mathfrak{M}\Gamma$.

We choose some $\omega$ satisfying
\begin{equation}\label{M1}
\beta+\gamma^{-\nu}\,a+\frac{1}{4}<\omega<\beta+\gamma^{-\nu}\,b-\frac{1}{4}
\end{equation}
and
\begin{equation}\label{M2}
\left|\frac{\omega}{2\pi}-\frac{q}{p}\right|\geq \frac{\gamma^\kappa}{p^\mu}
\end{equation}
for all integers $p,q$ with $p>0$ and $\mu>2$, $0<\kappa<\frac{1}{2}$. First we must show that for any sufficiently small $0<\gamma\leq 1$, there exists some $\omega$ satisfying (\ref{M1}) and (\ref{M2}). That is

\begin{lemma} For sufficiently small $0<\gamma\le 1$, there exists an $\omega$ satisfying (\ref{M1}) and (\ref{M2}).
\end{lemma}

\Proof Consider the complementary set $\Sigma$ of those $\omega$ in the interval $\Delta:=(\beta+\gamma^{-\nu}\,a+\frac{1}{4},\beta+\gamma^{-\nu}\,b-\frac{1}{4})$ that violate (\ref{M2}) for at least one pair of integers $p,q$ with $p\geq 1$. To estimate the Lebesgue measure of $\Sigma$, we fix $p$ and consider all $q$ for which the interval
$$
\left|\frac{\omega}{2\pi}-\frac{q}{p}\right|<\frac{\gamma^\kappa}{p^\mu}
$$
intersects the interval $\Delta$. Since the length of the interval $\Delta$ is $\gamma^{-\nu}(b-a)-\frac{1}{2}$, it is clear that there will be at most $\left[\frac{\gamma^{-\nu}(b-a)p}{2\pi}\right]+3$ such integers $q$, so that the measure of $\Sigma$ can be estimated by
$$
m(\Sigma)<\dsum_{p=1}^{+\infty} \left(\frac{\gamma^{-\nu}(b-a)p}{2\pi}+3\right)\frac{2\gamma^\kappa}{p^\mu}<8\gamma^{-\nu+\kappa}(b-a)\dsum_{p=1}^{+\infty} \frac{1}{p^{\mu-1}}=O(\gamma^\kappa)m(\Delta),
$$
which can be less than $m(\Delta)$ by taking $\gamma$ sufficiently small, here the assumption $b-a\ge 1$ is used. Hence the set $\Delta-\Sigma$ is not empty, and there exists an $\omega$ satisfying (\ref{M1}) and (\ref{M2}). The proof is completed. \qed

\medskip

Now we are in a position to  state our main result.

\begin{theorem}\label{A}
Under these hypotheses, for each sufficiently small $\gamma>0$ and any number $\alpha\in(\frac{1}{2},1)$, there exists a positive constant $c$, depending on $D$, $\alpha$ and $\mu$ in (\ref{M2}), but not on $\gamma$, $\nu$, $\kappa$ and $\beta$, such that for
\begin{equation}\label{small}
\gamma^{\nu+\kappa}|f|+|g|<c\,\gamma^{\frac{\nu+\kappa}{1-\alpha}}
\end{equation}
in $D$ the mapping $\mathfrak{M}$ in (\ref{M}) admits an invariant curve of the form
\begin{equation}\label{Curve}
x=\xi+ p(\xi;\gamma), \ \ \ \ \ \ \ \ y=q(\xi;\gamma)
\end{equation}
with $p,q$ real analytic functions of period $2\pi$ in the complex domain $|\Im\, \xi|<\frac{r_0}{2}$ and continuous in $\gamma$. Moreover, the parametrization is chosen so that the induced mapping on the curve (\ref{Curve}) is given by
\begin{equation}\label{RN}
\xi_1=\xi+\omega
\end{equation}
with $\omega$ satisfying (\ref{M1}) and (\ref{M2}), and the functions $p,q$ satisfy
$$
\gamma^{\nu+1}|p|+|q-\gamma^{\nu}\,(\omega-\beta)|<\bar{c}\,\gamma^{\nu+\frac{\nu+\kappa}{1-\alpha}},
$$
where $\bar{c}$ is same as $c$, depending on $D$, $\alpha$ and $\mu$ in (\ref{M2}), but not on $\gamma$, $\nu$, $\kappa$ and $\beta$.
\end{theorem}

Before turning to the proof of this theorem, it is useful to give some remarks on this result.

\begin{remark}
For $\nu = 0$, the mapping $\mathfrak{M}$ is the standard twist mapping, Theorem \ref{A} is the classical Moser twist theorem.
\end{remark}

\begin{remark} The unperturbed mapping of  $\mathfrak{M}$ is
\begin{equation}\label{UM}
x_1=x+\beta +\gamma^{-\nu}\, y, \quad\quad y_1=y,
\end{equation}
where the angle of rotation $\beta+\gamma^{-\nu}\, y$ ranges over the interval $[\beta+\gamma^{-\nu}\,a,\beta+\gamma^{-\nu}\,b]$, which will be large for small values of $\gamma$ and $\nu>0$. This is the reason we call it the large twist theorem.
\end{remark}

\begin{remark}\label{remark} The result is also valid for the general mapping
$$
\begin{array}{ll}
\left\{\begin{array}{ll}
x_1=x+\beta +\gamma^{-\nu}\,\alpha(y;\gamma)+f(x,y;\gamma),\\[0.2cm]
y_1=y+g(x,y;\gamma),\\[0.1cm]
 \end{array}\right.\  \ \ \ \ (x,y)\in \mathbb{R} \times [a,b],
\end{array}
$$
where $\alpha'(y;\gamma)\geq c_0$ for all $y\in [a,b]$ and the sufficiently small parameter $\gamma\in(0,1]$ with some constant $c_0\geq 0$. Indeed, if we introduce
$$x=\theta,\ \ \ \ \alpha(y;\gamma)=r
$$
as independent variables, then
$$
\alpha(a;\gamma)\leq r\leq \alpha(b;\gamma)
$$
and the mapping becomes
$$
\begin{array}{ll}
\left\{\begin{array}{ll}
\theta_1=\theta+\beta +\gamma^{-\nu}\, r+F(\theta,r;\gamma),\\[0.2cm]
r_1=r+G(\theta,r;\gamma),\\[0.1cm]
 \end{array}\right.\  \ \ \ \ (\theta,r)\in \mathbb{R} \times [\alpha(a;\gamma),\alpha(b;\gamma)]
\end{array}
$$
with
$$
F(\theta,r;\gamma)=f(x,y;\gamma),\ \ \ \ G(\theta,r;\gamma)=\alpha(y+g(x,y;\gamma);\gamma)-\alpha(y;\gamma).
$$
The width of the annulus can be estimated by
$$
\alpha(b;\gamma)-\alpha(a;\gamma)\geq c_0(b-a)\geq  c_0,
$$
and the smallness condition is same as (\ref{small}) with $c$ depending on $c_0$ more.
\end{remark}

\begin{remark}
In order to applying Theorem \ref{A} to the Poincar\'{e} mapping (\ref{P}), we need to verify that the functions $\Xi_1, \Xi_2$ in (\ref{P}) satisfy the smallness condition (\ref{small}) in Theorem \ref{A}. To this end, we choose $\alpha$ as
$$
\frac{1}{2}<\alpha<\frac{\nu+2\kappa}{2\nu+3\kappa}
$$
so that the smallness condition (\ref{small}) becomes
$$
|f|<c\,\gamma^{\nu+2\kappa}=c\,\gamma^{1-2\kappa}\,\gamma^{\nu+1}
$$
and
$$
|g|<c\,\gamma^{2\nu+3\kappa}=c\,\gamma^{\nu+3\kappa-1}\,\gamma^{\nu+1},
$$
where $\nu=n\geq 1$, and $1-2\kappa>0, \nu+3\kappa-1>0$.
\end{remark}

\begin{remark} As in the classical case (\cite{Moser}, \cite{Herman83}, \cite{Herman86}, \cite{Russmann83}), the analyticity of $f$ and $g$ can be replaced by some finite smoothness. For example, according to Herman (\cite{Herman83}, \cite{Herman86}),  R\"{u}ssman (\cite{Russmann83}), the smoothness requirement on  $f$ and $g$ can be weakened to $C^{3+\epsilon}$ for $\epsilon>0$.
\end{remark}

The proof of Theorem \ref{A} is based on the KAM approach, and to first find a sequence of coordinate changes
such that the transformed mapping of $\mathfrak{M}$ will be more closer to the unperturbed mapping (\ref{UM})
than the previous one in the narrower domain, then to make the infinite sequence of coordinate changes, with the range of these coordinates at the same time restricted to an annulus that shrinks down to the desired curve.

In the following we give a construction of such transformation. First, we choose some $\omega$ satisfying (\ref{M1}) and (\ref{M2}), consider the mapping $\mathfrak{M}$ in the complex domain
\begin{equation}\label{D}
D(r,s)\ :\ |\Im\,x|<r,\ \ \ |y-\bar{\omega}|<\gamma^{\nu}s
\end{equation}
with $\bar{\omega}=\gamma^{\nu}\,(\omega-\beta)$, $0<r<1$ and $0<s<{1\over 4}$, and construct a change of variables
\begin{equation}\label{UV}
\begin{array}{ll}
\mathfrak{U}:\quad \left\{\begin{array}{ll}
x = \xi+u(\xi,\eta),\\[0.2cm]
y = \eta+v(\xi,\eta),
\end{array}\right.
\end{array}
\end{equation}
where $u$ and $v$ are real analytic and  periodic functions in $\xi$.  Here and hereafter, the dependence of the functions on $\gamma$ is not indicated. Under this transformation, the original mapping $\mathfrak{M}$ is changed into the form
\begin{equation}\label{NewM}
\begin{array}{ll}
\mathfrak{{U}}^{-1} \mathfrak{M}\, \mathfrak{U}\ :\quad \ \left\{\begin{array}{ll}
\xi_1=\xi+\beta +\gamma^{-\nu}\,\eta+f_{+}(\xi,\eta),\\[0.2cm]
\eta_1=\eta+g_{+}(\xi,\eta),
 \end{array}\right.
\end{array}
\end{equation}
where the functions $f_{+}$ and $g_{+}$ are defined in a smaller domain $D(r_{+},s_{+})$ like (\ref{D}) and $\gamma^{\nu+\kappa}|f_{+}|_{r_+,s_+}+|g_{+}|_{r_+,s_+}$ is much smaller than $\gamma^{\nu+\kappa}|f|_{r,s}+|g|_{r,s}$, here $|f|_{r,s}=\sup\{|f(x,y)|: (x,y)\in D(r,s)\}$, and similar for the others.

From now on, $c_1,c_2,\cdots$ are positive constants depending on $\mu$ in (\ref{M2}) only.
From (\ref{M}), (\ref{UV}), (\ref{NewM}), it follows that
\begin{equation}\label{c4}
\left\{\begin{array}{ll}
f_{+}(\xi, \eta) =  f(\xi + u, \eta +  v) + \gamma^{-\nu}\,v(\xi, \eta) + u(\xi, \eta) - u(\xi_1,\eta_1),\\[0.3cm]
g_{+}(\xi, \eta) =  g(\xi + u, \eta +  v) + v(\xi, \eta) - v(\xi_1,\eta_1).
\end{array}\right.
\end{equation}
Now we will determine the unknown functions $u$ and $v$ from the following
equations
\begin{eqnarray*}
&&u(\xi + \omega, \eta)-u(\xi, \eta) = \gamma^{-\nu} \,v(\xi, \eta) +  f(\xi,\eta), \\[0.2cm]
&&v(\xi + \omega, \eta)-v(\xi, \eta) =  g(\xi, \eta).
\end{eqnarray*}
In order to obtain the analytic solutions of these difference equations, the problem of the small divisors is met. We can solve the second equation only if the mean value of $g$ over  the first variable vanishes. For this reason we define $u,v$ as the solutions of the following modified homological equations
\begin{equation}\label{c5}
\left\{\begin{array}{ll}
u(\xi + \omega, \eta)-u(\xi, \eta) = \gamma^{-\nu} \,v(\xi, \eta) +  f(\xi,\eta),\\[0.2cm]
v(\xi + \omega, \eta)-v(\xi, \eta) =   g(\xi, \eta)-  [g](\eta).
\end{array}\right.
\end{equation}
Here [\ ] denotes the mean value of a function over the first variable.\

We solve the functions $u$ and $v$ from (\ref{c5}) and give the estimates. In the first equation of (\ref{c5}), the  mean value over  the first variable must vanish on both sides. Hence we get the condition
\begin{equation*}\label{c6}
  [v](\eta)= -\gamma^{\nu}\,[f](\eta)
\end{equation*}
for  the  mean value of $v$ over  the first variable. As a consequence,  we have
\begin{equation*}
|[v]|_s\leq \gamma^{\nu}\, |f|_{r,s}.
\end{equation*}

\begin{lemma}\label{homogical equation}
Suppose $\omega$ satisfies (\ref{M1}) and (\ref{M2}). The difference equation
\begin{equation}\label{homological}
H(\xi+\omega,\eta)-H(\xi,\eta) = h(\xi,\eta)
\end{equation}
has a unique solution $H$ with $[H] = 0$, provided that the function $h$ is analytic in the domain $D(r,s)$ and $[h] = 0$. Moreover, the function $H$ is analytic in the domain $D(r-\delta,s)$ and the following estimate holds:
$$
|H|_{r-\delta,s}\leq c_1\gamma^{-\kappa} \delta^{-\mu}\,|h|_{r,s}.
$$
\end{lemma}

\Proof Suppose
$$
h(\xi,\eta) = \sum_{k\in\mathbb{Z}}h_k(\eta)e^{ik\xi}.
$$
Since the function $h(\xi,\eta)$ is analytic in $D(r,s)$ and $[h] = 0$, we have
$$
h_0(\eta) = 0,\quad |h_k|\leq |h|_{r,s}e^{-|k|r}.
$$
As usual, we have
$$
H(\xi,\eta) = \sum_{k\ne 0}\frac{h_k(\eta)}{e^{ik\omega}-1}e^{ik\xi}.
$$
Hence,
$$
\left|H(\xi,\eta)\right|_{r-\delta,s}\leq 2|h|_{r,s}\sum_{k\geq 1}\frac{1}{|e^{ik\omega-1}|}e^{-k\delta}.
$$
Since $\omega$ satisfies (\ref{M2}), it is easy to see that
$$
|e^{ik\omega}-1|\geq 4\frac{\gamma^{\kappa}}{k^{\mu-1}}.
$$
Hence we have
$$
\left|H(\xi,\eta)\right|_{r-\delta,s}\le  |h|_{r,s}\gamma^{-\kappa}\sum_{k=1}^\infty k^{\mu-1}e^{-k\delta}\le c_1\gamma^{-\kappa}\delta^{-\mu} \, |h|_{r,s}.
$$
The proof is finished.\qed

\medskip

By the above lemma, the second equation of (\ref{c5}) has a unique solution $\tilde{v}(\xi,\eta)$ with $[\tilde{v}](\eta)=0$, and this solution has the estimate
$$|\tilde{v}|_{r-\delta,s}\leq c_1\gamma^{-\kappa} \delta^{-\mu}\,|g|_{r,s}$$
for $0<\delta<r.$ Define $v(\xi,\eta)=\tilde{v}(\xi,\eta)+[v](\eta)$, we obtain the  uniquely determined solution $v(\xi,\eta)$ of the second equation of (\ref{c5}).

Define $p(\xi,\eta)=\gamma^{-\nu} \,\tilde{v}(\xi,\eta)+ f(\xi,\eta)$, note that  $\tilde{v}(\xi,\eta)$ is defined in $D(r-\delta,s)$,\ then  $p(\xi,\eta)$ is well defined in $D(r-\delta,s)$. As a consequence we have
\begin{eqnarray*}
|p|_{r-\delta,s} =|\gamma^{-\nu}\,\tilde{v}+ f|_{r-\delta,s}
\leq  c_2 \gamma^{-\nu-\kappa}\delta^{-\mu}\,(\gamma^{\nu+\kappa}|f|_{r,s}+|g|_{r,s})
\end{eqnarray*}
and
\begin{eqnarray*}
p(\xi,\eta)-[p](\eta)
&=&\gamma^{-\nu}\, \tilde{v}(\xi,\eta)+f(\xi,\eta)-\gamma^{-\nu}\,[ \tilde{v}](\eta)-  [f](\eta) \\[0.2cm]
&=&  \gamma^{-\nu}\, \tilde{v}(\xi,\eta)+ f(\xi,\eta)+\gamma^{-\nu}\, [v](\eta)\\[0.2cm]
&= & \gamma^{-\nu}\,  v(\xi,\eta)+ f(\xi,\eta).
\end{eqnarray*}
Hence, the first equation of (\ref{c5}) can  be written in  the  form
\begin{equation}\label{c7}
u(\xi + \omega, \eta)-u(\xi, \eta) = p(\xi, \eta) - [p](\eta).
\end{equation}
Using Lemma \ref{homogical equation} again,  (\ref{c7}) has a uniquely determined solution  $u$  with $[u]=0,$ which possesses the estimate
$$|u|_{r-2\delta,s}\leq c_1\gamma^{-\kappa}\delta^{-\mu}|p|_{r-\delta,s}\leq  c_1c_2\gamma^{-\nu-2\kappa} \delta^{-2\mu}\,(\gamma^{\nu+\kappa}|f|_{r,s}+|g|_{r,s}).$$

From the above discussions, we have
\begin{eqnarray}\label{c8}
\gamma^{\nu+\kappa}|u|_{r-2\delta,s}+|v|_{r-2\delta,s} &\leq& c_{3}\gamma^{-\kappa}\delta^{-2\mu}\,(\gamma^{\nu+\kappa}|f|_{r,s}+|g|_{r,s})\nonumber\\[0.2cm]
&:=& c_{3}\gamma^{-\kappa} \delta^{-2\mu} d
\end{eqnarray}
with
$$d=\gamma^{\nu+\kappa}|f|_{r,s}+|g|_{r,s},$$
and by Cauchy's estimate,
\begin{equation}\label{c70}
\gamma^{\nu+\kappa}\left|D_{\xi} u\right|_{r-3\delta,s}+\left|D_{\xi} v\right|_{r-3\delta,s} \leq  c_3\gamma^{-\kappa} \delta^{-2\mu-1} d,
\end{equation}
\begin{equation}\label{c71}
\gamma^{\nu+\kappa}\left|D_{\eta} u\right|_{r-2\delta,s-\rho}+\left|D_{\eta} v\right|_{r-2\delta,s-\rho}
 \leq  c_3\gamma^{-\kappa} \delta^{-2\mu} {d\over \rho}
\end{equation}
for $0<3\delta<r$ and $0<\rho<s$.

For $0<r_+<r$ and $0<s_+<s$, let
$$\delta={1\over 6}(r-r_+),\ \ \ \ \rho={1\over 6}(s-s_+).$$
Introduce intermediate domains $D_j$ between $D(r_+,s_+)$ and $D(r,s)$ by
$$D_j=D(r-j\delta,s-j\rho)\ \ \ \ \mbox{for}\ \ 0\leq j\leq 6.$$
From (\ref{c8}, (\ref{c70}) and (\ref{c71}), it follows that
\begin{equation*}\label{c80}
\gamma^{\nu+\kappa}|u|_{D_2}+|v|_{D_2} \leq c_{4}\gamma^{-\kappa}(r-r_+)^{-2\mu}\,d,
\end{equation*}
\begin{equation*}\label{c81}
\gamma^{\nu+\kappa}\left|D_{\xi} u\right|_{D_3}+\left|D_{\xi} v\right|_{D_3} \leq  c_4\gamma^{-\kappa} (r-r_+)^{-2\mu-1} d,
\end{equation*}
and
\begin{equation*}\label{c82}
\gamma^{\nu+\kappa}\left|D_{\eta} u \right|_{D_2}+\left|D_{\eta} v\right|_{D_2}
 \leq  c_4\gamma^{-\nu-\kappa}  (r-r_+)^{-2\mu}\, {d\over s},
\end{equation*}
where $c_4>c_3$, and we also need that $7s_+<s.$ With this constant $c_4$ we define
$$
\vartheta=c_4 \gamma^{-\nu-\kappa}(r-r_+)^{-2\mu}\, {d\over s}
$$
and
rewrite these inequalities in the form
\begin{equation}\label{DofUV}
\left\{\begin{array}{ll}
\gamma^{\nu+\kappa}|u|_{D_2}+|v|_{D_2} \leq \gamma^{\nu}\vartheta s,\\[0.4cm]
\gamma^{\nu+\kappa}\left|D_{\xi} u\right|_{D_3}+\left|D_{\xi} v\right|_{D_3}  \leq \gamma^{\nu}\vartheta \frac{s}{r-r_+}<\gamma^{\nu}\vartheta,\\[0.4cm]
\gamma^{\nu+\kappa}\left|D_{\eta} u \right|_{D_2}+\left|D_{\eta} v\right|_{D_2}
 \leq \vartheta,
\end{array}\right.
\end{equation}
where we use that
$$
s<r-r_+.
$$

Therefore, since the parameter $\gamma\in(0,1]$, similar to the classical case (\cite{Siegel97}), if
$$
0<7s_+<s<\frac{\gamma^\kappa}{6}(r-r_+),\ \ \ \ d<\frac{\gamma^{\nu+\kappa}}{7}s, \ \ \ \ \vartheta<\frac{1}{8},
$$
one can prove that
\begin{eqnarray*}\label{c10}
\mathfrak{U}(D_6)\subset D_5,\ \ \ \ \mathfrak{M}(D_{5})\subset D_4,\ \ \ \ \mathfrak{U}^{-1}(D_{4})\subset D_3.
\end{eqnarray*}

In fact, in order to obtain $\mathfrak{U}(D_6)\subset D_5$, we need to verify that
$$
(|u|_{D_2}  \leq) \gamma^{-\kappa}\vartheta s\leq \frac{r-r_+}{6}
$$
and
$$
(|v|_{D_2} \leq) \gamma^{\nu}\vartheta s \leq\frac{\gamma^{\nu}}{6}(s-s_+),
$$
which are guaranteed by $7s_+<s<\frac{\gamma^\kappa}{6}(r-r_+), \vartheta<\frac{1}{8}$. Similarly, to check that $\mathfrak{M}$ maps $D_{5}$ into $D_{4}$,  we use the expression of $\mathfrak{M}$ and find that
$$
\begin{array}{lll}
x_1&=&x+\beta +\gamma^{-\nu}\,y+ f\\[0.2cm]
   &=&x+\beta +\gamma^{-\nu}\bar{\omega}+\gamma^{-\nu}\,(y-\bar{\omega})+ f\\[0.2cm]
 &=&x+\omega+\gamma^{-\nu}\,(y-\bar{\omega})+ f,\\[0.2cm]
y_1&=&y+g,
\end{array}
$$
which implies that
$$
|\Im x_1|\leq r-5\delta+s-5\rho+\gamma^{-\nu-\kappa}d,
$$
and
$$
|y_1-\bar{\omega}|<|y-\bar{\omega}|+d<\gamma^{\nu}(s-5\rho)+d,
$$
so that we only have to verify the inequalities
$$
s-5\rho+\gamma^{-\nu-\kappa}d<\frac{r-r_+}{6},\ \ \ \  d<\frac{\gamma^{\nu}}{6}(s-s_+),
$$
which are also guaranteed by $7s_+<s<\frac{\gamma^\kappa}{6}(r-r_+), d<\frac{\gamma^{\nu+\kappa}}{7}s$.

Finally, we show that $\mathfrak{U}^{-1}$ is defined in the domain $D_{4}$ and maps it into $D_3$. In other word, we assume that $(x,y)\in D_4$ and have to construct a solution $(\xi, \eta)$ of the equation (\ref{UV}) in $D_3$. For this purpose we use the usual iteration scheme and define $\xi_k, \eta_k$ inductively by $\xi_0=x, \eta_0=y$ and
$$
\xi_{k+1}=x-u(\xi_k, \eta_k),\ \ \eta_{k+1}=y-v(\xi_k, \eta_k), \ \ k=0,1,\cdots.
$$
We must show that the $(\xi_k, \eta_k)$ remain in $D_3$, that is, $|\Im \xi_k|<r-3\delta, |\eta_k-\bar{\omega}|<\gamma^{\nu}(s-3\rho)$. For $k=0$ this is obviously the case, and assume this to be true for $(\xi_l, \eta_l)$ with $l\leq k$, we find from (\ref{DofUV}) that
$$
\begin{array}{lll}
\gamma^{\nu+\kappa}|\xi_{k+1}-\xi_k|+|\eta_{k+1}-\eta_k|&\leq& \left(\gamma^{\nu+\kappa}\left|D_{\xi} u \right|_{D_3}+\left|D_{\xi} v\right|_{D_3}\right)|\xi_{k}-\xi_{k-1}|\\[0.4cm]
&+& \left(\gamma^{\nu+\kappa}\left|D_{\eta} u \right|_{D_3}+\left|D_{\eta} v\right|_{D_3}\right)|\eta_{k}-\eta_{k-1}|\\[0.4cm]
&\leq& \vartheta \left(\gamma^{\nu}|\xi_{k}-\xi_{k-1}|+|\eta_{k}-\eta_{k-1}|\right)\\[0.4cm]
&\leq& \vartheta \left(\gamma^{\nu+\kappa}|\xi_{k}-\xi_{k-1}|+|\eta_{k}-\eta_{k-1}|\right)
\end{array}
$$
and hence
$$
\begin{array}{lll}
\gamma^{\nu+\kappa}|\xi_{k+1}-x|+|\eta_{k+1}-y|&\leq& \dsum_{l=1}^{k+1} \left(\gamma^{\nu+\kappa}|\xi_{l}-\xi_{l-1}|+|\eta_{l}-\eta_{l-1}|\right)\\[0.4cm]
&\leq& \frac{1}{1-\vartheta}\left(\gamma^{\nu+\kappa}|u|_{D_3}+|v|_{D_3}\right)\\[0.4cm]
&\leq& \frac{\vartheta}{1-\vartheta}\gamma^{\nu}s.
\end{array}
$$
Since we assume that $\vartheta<\frac{1}{8}$, the last inequality is $<\frac{\gamma^{\nu}}{7}s$ and as before, we verify
that
$$
\frac{\gamma^{-\kappa}}{7}s<\frac{r-r_+}{6},\ \ \ \ \frac{\gamma^{\nu}}{7}s<\frac{\gamma^{\nu}}{6}(s-s_+),
$$
which guarantees that  $(\xi_{k+1}, \eta_{k+1})\in D_3$. Thus all the iterates $(\xi_{k}, \eta_{k})$ remain in $D_3$ and since $\vartheta<1$, converge to a solution $(\xi, \eta)$ in $D_3$. This solution is clearly unique, and $\mathfrak{U}^{-1}$ maps $D_{4}$ into $D_3$ as asserted.

Hence the mapping $\mathfrak{U}^{-1}\mathfrak{M}\,\mathfrak{U}$ is well defined in the domain $D(r_+,s_+)$ and maps this domain into $D_3$.
Similar to the proof in \cite[chapter 3]{Siegel97}, one may prove that $f_+$ and $g_+$ defined in (\ref{c4}) are real analytic in $D(r_+,s_+)$ and for each fixed $\eta$, periodic in $\xi$.

In the following, we will prove that $\gamma^{\nu+\kappa}|f_+|_{r_+,s_+}+|g_+|_{r_+,s_+}$ is much smaller than $\gamma^{\nu+\kappa}|f|_{r,s}+|g|_{r,s}$.

From (\ref{c4}) and (\ref{c5}), we have
\begin{equation*}
\begin{array}{ll}
f_{+}(\xi, \eta) = u(\xi+\omega, \eta) - u(\xi_1,\eta_1)+f(\xi + u, \eta +v) -f(\xi,\eta) ,\\[0.3cm]
g_{+}(\xi, \eta) = v(\xi+\omega, \eta)- v(\xi_1,\eta_1)+g(\xi + u, \eta + v)-g(\xi,\eta)+[g](\eta).
\end{array}
\end{equation*}

The troublesome mean value  $[g](\eta)$ will be approximated by the linear function
$$h(\eta)=[g](\bar{\omega})+D_{\eta}[g](\bar{\omega})(\eta-\bar{\omega}),$$
which we will estimate later using the intersection property. A preliminary estimate of $h(\eta)$ and $\left|[g]-h\right|_{s_+}$ is obtained by observing that for $|\eta-\bar{\omega}|<\gamma^{\nu}s$ we have $\left|[g]\right|_{s}<d$ and therefore by Cauchy's estimate $\left|D_{\eta}[g](\bar{\omega})\right|<\frac{d}{\gamma^{\nu}s},$ while for $|\eta-\bar{\omega}|<\gamma^{\nu}s_+$ also $\left|D_{\eta}^2[g]\right|_{s_+}<{{2d}\over{\gamma^{2\nu}(s-s_+)^2}}.$ Consequently,  for $|\eta-\bar{\omega}|<\gamma^{\nu}s_+$, we have
$$|h|_{s_+}<d+{d\over{s}}s_+<2d$$
and
$$\left|[g]-h\right|_{s_+}<{{\gamma^{2\nu}s_+^2}\over 2}\left|D_{\eta}^2[g]\right|_{s_+}<\left({{s_+}\over{s-s_+}}\right)^2 d<2\left(\frac{s_+}{s}\right)^2d,$$
where we have used that $7s_+<s$.

By Cauchy's estimate and (\ref{DofUV}), it follows that
$$
\begin{array}{lll}
&&|g_+ -h|_{r_+,s_+}\\[0.2cm]
&\leq & |v(\xi+\omega,\eta)-v(\xi_{1},\eta_{1})|_{r_+,s_+}+|g(\xi+u,\eta+v)-g(\xi,\eta)|_{r_+,s_+}+|[g]-h|_{s_+} \\[0.2cm]
&\leq & \gamma^{\nu}\vartheta \frac{s^2}{r-r_+}+\vartheta (\gamma^{\nu+\kappa}|f_+|_{r_+,s_+}+|g_+|_{r_+,s_+})+2\vartheta d +2\left(\frac{s_+}{s}\right)^2d\\[0.2cm]
&\leq& \gamma^{\nu}\vartheta\frac{s^2}{r-r_+}+\vartheta \left(\gamma^{\nu+\kappa} |f_+|_{r_+,s_+}+|g_+ -h|_{r_+,s_+}+|h|_{s_+}\right)+2 \vartheta d +2\left(\frac{s_+}{s}\right)^2d \\[0.2cm]
&\leq& \gamma^{\nu}\vartheta\frac{s^2}{r-r_+}+\vartheta \left( \gamma^{\nu+\kappa}|f_+|_{r_+,s_+}+|g_+ -h|_{r_+,s_+}+2d\right)+ 2\vartheta d +2\left(\frac{s_+}{s}\right)^2d.
\end{array}
$$
For $f_{+}$, similarly we have
$$
\begin{array}{lll}
&&\gamma^{\nu+\kappa} |f_+|_{r_+,s_+}\\[0.2cm]
&\leq& \gamma^{\nu+\kappa} |u(\xi+\omega,\eta)-u(\xi_{1},\eta_{1})|_{r_+,s_+}+\gamma^{\nu+\kappa} |f(\xi+u,\eta+ v)-f(\xi,\eta)|_{r_+,s_+}\\[0.2cm]
&\leq & \gamma^{\nu}\vartheta \frac{s^2}{r-r_+}+\vartheta (\gamma^{\nu+\kappa} |f_+|_{r_+,s_+}+|g_+|_{r_+,s_+})+2\vartheta d \\[0.2cm]
&\leq& \gamma^{\nu}\vartheta\frac{s^2}{r-r_+}+\vartheta \left(\gamma^{\nu+\kappa}  |f_+|_{r_+,s_+}+|g_+ -h|_{r_+,s_+}+2d\right)+ 2\vartheta d .
\end{array}
$$
Adding this two estimates we obtain
$$
\begin{array}{lll}
& \gamma^{\nu+\kappa} |f_+|_{r_+,s_+}+|g_+ -h|_{r_+,s_+}\\[0.2cm]
\leq& 2\gamma^{\nu}\vartheta\frac{s^2}{r-r_+}+2\vartheta \left(\gamma^{\nu+\kappa}  |f_+|_{r_+,s_+}+|g_+ -h|_{r_+,s_+}+2d\right)+ 4\vartheta d +2\left(\frac{s_+}{s}\right)^2d.
\end{array}
$$
If we assume that $\vartheta<\frac{1}{8}$, then we  can eliminate $\gamma^{\nu+\kappa} |f_+|_{r_+,s_+}+|g_+ -h|_{r_+,s_+}$ from the right hand side, and obtain
\begin{align}\label{c17}
 \gamma^{\nu+\kappa} |f_+|_{r_+,s_+}+|g_+ -h|_{r_+,s_+} < c_5\left(\frac{\gamma^{\nu}\vartheta }{r-r_+}(s^2+\gamma^{-\nu}d)+\left(\frac{s_+}{s}\right)^2d\right),
\end{align}
which implies together with the express of $\vartheta$ that
\begin{align}\label{c171}
 \gamma^{\nu+\kappa} |f_+|_{r_+,s_+}+|g_+ -h|_{r_+,s_+} < Q,
\end{align}
where
\begin{align}\label{c1700}
Q:=c_5\left\{c_4\gamma^{-\kappa}(r-r_+)^{-2\mu-1}\left(sd+\gamma^{-\nu}\frac{d^2}{s}\right) +\left(\frac{s_+}{s}\right)^2d\right\}.
\end{align}

The preliminary estimate of $2d$ for $|h|_{s_+}$, however, is insufficient for decreasing the error term, and to obtain a better estimate we use the
intersection property of $\mathfrak{M}$, or that of $\mathfrak{N}=\mathfrak{{U}}^{-1} \mathfrak{M}\, \mathfrak{U}$. Accordingly, each curve  $\eta=\text{constant}$, in particular, has to intersect its image curve under $\mathfrak{N}$ and at such a point of intersection we have $\eta_1=\eta$ or $g_+=0$, so that for each real $\eta$ in $|\eta-\bar{\omega}|<\gamma^{\nu}s_+$, there exists a real $\xi=\xi_0(\eta)$ such that $g_+(\xi_0(\eta),\eta)=0.$ Applying (\ref{c171}) at such points $(\xi_0(\eta),\eta)$, we find that
$$|h(\eta)|<Q\ \ \ (\bar{\omega}-\gamma^{\nu}s_+<\eta<\bar{\omega}+\gamma^{\nu}s_+).$$
Consequently, setting $\eta=\bar{\omega}$ in the definition of $h$, we get $\big|[g](\bar{\omega})\big|<Q,$ and
letting  $\eta$ approach $\bar{\omega}+\gamma^{\nu}s_+$ in the same definition, we obtain
$$\big|[g](\bar{\omega})+D_{\eta}[g](\bar{\omega})\gamma^{\nu}s_+\big|\leq Q,$$
so that
$$\big|D_{\eta}[g](\bar{\omega})\gamma^{\nu}s_+\big|< 2Q.$$
From this we conclude that for complex $\eta$ in the disk $|\eta-\bar{\omega}|<\gamma^{\nu}s_+$ we have
$$\big|h(\eta)\big|\leq \big|[g](\bar{\omega})\big|+\big|D_{\eta}[g](\bar{\omega})\big|\big|\eta-\bar{\omega}\big|<3Q,$$
which in view of (\ref{c171}) gives
$$\gamma^{\nu+1} |f_+|_{r_+,s_+}+|g_+|_{r_+,s_+}<4Q.$$

The above discussions lead to the following lemma.

\begin{lemma}\label{lem3.1}
Consider the map $\mathfrak{M}$ defined by (\ref{M}), where $f$ and $g$ are real analytic in the domain $D(r,s)$ and periodic in $x$, and the frequency $\omega$ satisfies (\ref{M1}) and (\ref{M2}). Denote
$$d=\gamma^{\nu+\kappa} |f|_{r,s}+|g|_{r,s}.$$
Assume that the conditions
\begin{align}\label{c170}
0<7s_+<s<\frac{\gamma^\kappa}{6}(r-r_+),\ \ \ \ d<\frac{\gamma^{\nu+\kappa}}{7}s, \ \ \ \ \vartheta<\frac{1}{8}
\end{align}
are satisfied, then there exists a transformation
\begin{equation*}
\begin{array}{ll}
\mathfrak{U}:\ \ \left\{\begin{array}{ll}
x = \xi+u(\xi,\eta),\\[0.2cm]
y = \eta+v(\xi,\eta),
\end{array}\right.
\end{array}
\end{equation*}
which is defined in a smaller domain $D(r_+,s_+)$. Under this transformation, the map $\mathfrak{M}$ has the form
\begin{equation*}
\begin{array}{ll}
\mathfrak{{U}}^{-1} \mathfrak{M}\, \mathfrak{U}\ :\quad \ \left\{\begin{array}{ll}
\xi_1=\xi+\beta +\gamma^{-\nu}\eta+f_{+}(\xi,\eta),\\[0.2cm]
\eta_1=\eta+ g_{+}(\xi,\eta),
 \end{array}\right.
\end{array}
\end{equation*}
where the functions $f_{+}$ and $g_{+}$ are real analytic in the domain $D(r_{+},s_{+})$ and periodic in $\eta$. Moreover, the following estimates hold :
$$
\gamma^{\nu+\kappa} |u|_{D_2}+|v|_{D_2} \leq \gamma^{\nu}\vartheta s, \ \ \ \ \gamma^{\nu+\kappa}|f_+|_{r_+,s_+}+|g_+|_{r_+,s_+} < 4Q.
$$
\end{lemma}

Lemma \ref{lem3.1} is usually called the iteration lemma in the KAM proof. For the proof of Theorem \ref{A}, one can use this lemma infinite times to construct a sequence of transformations. To this end we make a sequence of successive applications of the lemma, starting with the given mapping  $\mathfrak{M}$ in (\ref{M}), now denoted by $\mathfrak{M}_0=\mathfrak{M}$ and restricted to the domain
$$
D_0\ :\ |\Im\,x|<r_0,\ \ \ |y-\bar{\omega}|<\gamma^{\nu}s_0,
$$
which for $s_0$ sufficiently small is contained in the domain $D$ given by (\ref{D1}), also we need to assume that $s_0<\frac{1}{4}$ so that whenever $\bar{\omega}\in (a+\frac{1}{4}\gamma^{\nu}, b-\frac{1}{4}\gamma^{\nu})$ determined by (\ref{M1}), any real $y$ satisfying $|y-\bar{\omega}|<\gamma^{\nu}s_0$ is contained in the interval $[a,b]$. We first assume that
$$
\gamma^{\nu+\kappa}|f|_{D_0}+|g|_{D_0}<d_0,
$$
the key problem is to find the relation between $d_0$ and $\gamma$. Transforming the mapping $\mathfrak{M}$ by the coordinate $\mathfrak{U}=\mathfrak{U}_0$ provided by the lemma, we obtain a mapping $\mathfrak{M}_1=\mathfrak{U}_0^{-1}\mathfrak{M}_0\mathfrak{U}_0$ defined in the domain
$$
D_1\ :\ |\Im\,x|<r_1,\ \ \ |y-\bar{\omega}|<\gamma^{\nu}s_1,
$$
where $r_1, s_1$ correspond to the parameter $r_+, s_+$ of the lemma. Applying the lemma to the new mapping $\mathfrak{M}_1$, we obtain another coordinate transformation $\mathfrak{U}_1$ and a transformed mapping $\mathfrak{M}_2=\mathfrak{U}_1^{-1}\mathfrak{M}_1\mathfrak{U}_1$, and proceeding in this way we are led to a sequence of mappings
\begin{equation}\label{Mn}
\mathfrak{M}_{n+1}=\mathfrak{U}_n^{-1}\mathfrak{M}_n\mathfrak{U}_n,\quad\quad n=0,1,\cdots,
\end{equation}
whose domains $D_{n+1}$ are defined like $D$ by (\ref{D}) with $r_{n+1}, s_{n+1}$ replacing $r,s$. We have to verify, of course, that this sequence of transformations is well defined and that $\mathfrak{M}_n$ approximates the unperturbed twist mapping with increasing precision. For this we fix the parameters $r_n, s_n, d_n$ $(n=0,1,\cdots)$ by setting
\begin{equation}\label{Mn}
r_n=\frac{r_0}{2}\left(1+\frac{1}{2^n}\right),\ \ \ \  s_n=d_n^{\alpha},\ \ \ \ d_{n+1}= \gamma^{-\nu-\kappa}r_0^{-2\mu-1} c_7^{n+1}d_n^{2-\alpha}
\end{equation}
with $c_7\geq 7$ a suitably chosen constant, $\frac{1}{2}<\alpha<1$ is any number. Thus $r_n$ is a decreasing sequence converging to the positive value $\frac{r_0}{2}$, and all functions to be considered will be analytic in $\xi$ for $|\Im \xi|<\frac{r_0}{2}$. The sequence $d_n$ converges to $0$ provided $d_0$ is chosen sufficiently small. Indeed the sequence
$$
e_n=\gamma^{-\frac{\nu+\kappa}{1-\alpha}}r_0^{-\frac{2\mu+1}{1-\alpha}}c_7^{\frac{n}{1-\alpha}+\frac{2-\alpha}{(1-\alpha)^2}}d_n
$$
satisfies
$$
e_{n+1}=e_n^{2-\alpha}
$$
and therefore converges to zero if we take $0\leq e_0<1$, or
\begin{equation}\label{d_0}
0\leq d_0<\gamma^{\frac{\nu+\kappa}{1-\alpha}}r_0^{\frac{2\mu+1}{1-\alpha}}c_7^{-\frac{2-\alpha}{(1-\alpha)^2}}.
\end{equation}
Moreover, $e_n=e_0^{(2-\alpha)^n}$, or
\begin{equation}\label{d_n}
d_n=\gamma^{\frac{\nu+\kappa}{1-\alpha}}r_0^{\frac{2\mu+1}{1-\alpha}}c_7^{-\frac{n}{1-\alpha}-\frac{2-\alpha}{(1-\alpha)^2}}e_0^{(2-\alpha)^n}, \ \ n=0,1,\cdots,
\end{equation}
and
\begin{equation}\label{s_n}
s_n=\gamma^{\frac{(\nu+\kappa)\alpha}{1-\alpha}}r_0^{\frac{(2\mu+1)\alpha}{1-\alpha}}c_7^{-\frac{n\alpha}{1-\alpha}-\frac{(2-\alpha)\alpha}{(1-\alpha)^2}}e_0^{\alpha(2-\alpha)^n}, \ \ n=0,1,\cdots
\end{equation}

To show that the mapping $\mathfrak{M}_n$ is well defined in $D_n$ and satisfies there the appropriate estimate, we proceed by induction. For $n=0$, we choose $r_0,s_0,d_0$ in (\ref{Mn}), (\ref{d_n}),(\ref{s_n}) with $e_0$ sufficiently small, it is easy to see that the condition (\ref{c170}) holds. Assuming $\mathfrak{M}_n$ to be defined in $D_n$ and satisfying the estimate
$$
\gamma^{\nu+\kappa}|f|_{D_n}+|g|_{D_n}<d_n,
$$
we will verify the corresponding estimate for $\mathfrak{M}_{n+1}$. For this we apply the lemma with $r=r_n,s=s_n,d=d_n,r_+=r_{n+1},s_+=s_{n+1}$, and therefore have to verify first the inequalities (\ref{c170}).  By (\ref{s_n}), it follows that
\begin{equation}\label{s_n1}
\left(\frac{s_{n+1}}{s_n}\right)^{\frac{1}{\alpha}}=\frac{d_{n+1}}{d_n}=c_7^{-\frac{1}{1-\alpha}}e_n^{1-\alpha}<c_7^{-2}\leq \frac{1}{49},
\end{equation}
which implies that $7s_+<s$. The inequality $s<\frac{\gamma^\kappa}{6}(r-r_+)$ follows from (\ref{s_n}) and
$$
r-r_+=r_n-r_{n+1}=r_02^{-n-2}
$$
with $e_0$ sufficiently small. Clearly also the inequality $d<\frac{\gamma^{\nu+\kappa}}{7}s$ can be met by a suitable restriction on $e_0$, while for $\vartheta=\vartheta_n$ we have
$$
\vartheta_n=c_4 \gamma^{-\nu-\kappa}(r_n-r_{n+1})^{-2\mu}\, d_n^{1-\alpha}
$$
and  this too can be made less than $\frac{1}{8}$ by choosing $e_0$ small. Thus there exists a positive constant $c=c(r_0,\mu,\alpha)$ such that for
\begin{equation}\label{d_01}
d_0<c\,\gamma^{\frac{\nu+\kappa}{1-\alpha}}
\end{equation}
the inequalities (\ref{c170}) hold and the lemma is applicable. From the lemma we now obtain the transformation $\mathfrak{U}_n$ taking $D_{n+1}$ into $D_n$ and the transformed mapping $\mathfrak{M}_{n+1}=\mathfrak{U}_n^{-1}\mathfrak{M}_n\mathfrak{U}_n$ defined in $D_{n+1}$. Moreover, representing $\mathfrak{M}_{n+1}$ in the form (\ref{NewM}), by the lemma we have the estimate
$$
\begin{array}{lll}
&& \gamma^{\nu+\kappa} |f_+|+|g_+|\\[0.1cm]
&<& 4c_5\left\{c_4\gamma^{-\kappa}(r_n-r_{n+1})^{-2\mu-1}\left(s_nd_n+\gamma^{-\nu}\frac{d_n^2}{s_n}\right) +\left(\frac{s_{n+1}}{s_n}\right)^2d_n\right\} \\[0.4cm]
&<& c_6\left\{\gamma^{-\nu-\kappa}(r_n-r_{n+1})^{-2\mu-1}\left(s_nd_n+\frac{d_n^2}{s_n}\right) +\left(\frac{s_{n+1}}{s_n}\right)^2d_n\right\} \\[0.4cm]
&=& c_6\left\{\gamma^{-\nu-\kappa}(r_n-r_{n+1})^{-2\mu-1}\left(d_n^{1+\alpha}+d_n^{2-\alpha}\right) +\left(\frac{d_{n+1}}{d_n}\right)^{2\alpha}d_n\right\} \\[0.4cm]
&=& c_6\left\{\gamma^{-\nu-\kappa}(r_n-r_{n+1})^{-2\mu-1}d_n^{2-\alpha}\left(1+d_n^{2\alpha-1}\right) +\left(\frac{d_{n+1}}{d_n}\right)^{2\alpha-1}d_{n+1}\right\} \\[0.4cm]
&=& c_6\left\{2^{(n+2)(2\mu+1)}c_7^{-(n+1)}\left(1+d_n^{2\alpha-1}\right)+\left(\frac{d_{n+1}}{d_n}\right)^{2\alpha-1}\right\}d_{n+1} \\[0.4cm]
&=& c_6\left\{2^{2\mu+1}\left(2^{2\mu+1}c_7^{-1}\right)^{n+1}\left(1+d_n^{2\alpha-1}\right)+c_7^{-\frac{2\alpha-1}{1-\alpha}}\right\}d_{n+1},
\end{array}
$$
and since the sequence $d_n$ is bounded, $\frac{1}{2}<\alpha<1$, the coefficient of $d_{n+1}$ in the last inequality can be made less than $1$ by taking $c_7$ large. With such a choice we finally have
$$
\gamma^{\nu+\kappa} |f_+|+|g_+|<d_{n+1},
$$
and the induction is complete.

The remaining part is completely same as the classical case (see \cite{{Siegel97}}), and we omit it here.

\section{The proof of Theorem \ref{B}}

In this section we first introduce the action-angle variables. Let $(C(t),S(t))$ be the solution of
$$
\dot{x}=y,\quad\quad \dot{y} =-x^{2n+1}
$$
with the initial value $(C(0),S(0)) = (1,0)$, and $T^*$ be its minimal period.

The action and angle variables are defined by the map $\Psi: \mathbb{R^+}\times \mathbb{S}\to \mathbb{R}^2\backslash\{0\}$
$$
\Psi: x=a^{\alpha}\lambda^{\alpha}C(\theta T^*), \quad\quad y=a^{\beta}\lambda^{\beta}S(\theta T^*)
$$
with $\alpha=\frac{1}{n+2}, \beta=1-\frac{1}{n+2}$, and $a=\frac{1}{\alpha T^*}$.

Introduce $y=\dot{x}$, then Eq.(\ref{Duffing}) is a Hamilton system with the Hamiltonian
$$
H(x,y,t)=\frac{1}{2}y^2+\frac{1}{2n+2}x^{2n+2}+\sum_{i=0}^{2n}\frac{p_i(t)}{i+1}x^{i+1},
$$
and under the action and angle variables, which has the form
\begin{equation}\label{H1}
h(\lambda,\theta,t) = d\cdot\lambda^{2\beta} + h_1(\lambda,\theta,t) + h_2(\lambda,\theta,t),
\end{equation}
where
$$
h_1(\lambda,\theta,t)=\sum_{i=n+1}^{2n}a_i\lambda^{(i+1)\alpha}C(\theta T^*)^{i+1}p_i(t),
$$
$$
h_2(\lambda,\theta,t)=\sum_{i=0}^{n}a_i\lambda^{(i+1)\alpha}C(\theta T^*)^{i+1}p_i(t),
$$
and $d = \frac{a^{2\beta}}{2n+2}, a_i=\frac{a^{(i+1)\alpha}}{i+1}$.

In order to perform some canonical transformations for the Hamiltonian $h$ defined by (\ref{H1}), we first introduce the function space similar to \cite{Diecherhoff}. For $r\in \mathbb{R}$ and $k\in \mathbb{N}$, we say a function $f(\lambda,\theta,t)\in F^k(r)$ if $f$ is smooth in $\lambda$ and $\theta$, $C^k$ in $t$, and
$$
\sup_{\lambda\ge\lambda_0,(\theta,t)\in \mathbb{S}\times \mathbb{S}}\lambda^{j-r}\left|D^j_\lambda D^l_\theta f(\lambda,\theta,t)\right| < +\infty.
$$
We summarize some properties readily verified from the definition.

\begin{lemma}
(1) if $r_1 < r_2$, then $F^k(r_1) \subset F^k(r_2)$;\\
(2) if $f\in F^k(r)$, then $D_{\lambda}^j f\in F^k(r - j)$;\\
(3) if $f_1\in F^k(r_1), f_2 \in F^k(r_2)$, then $f_1\cdot f_2\in F^k(r_1 + r_2)$;\\
(4) if $f\in F^k(r)$, and satisfies $|f(\lambda,\cdot)| \geq c\cdot\lambda^r$  for $\lambda\geq\lambda_0$, then $\frac{1}{f}\in F^k(-r)$;\\
(5) if $f\in F^k(r), g\in F^k(s)$, then $f(g)\in F^k(rs)$.
\end{lemma}

For $f\in F^k(r)$, we denote the mean value of $f$ over the $\theta$-variable by $[f]: [f](\lambda,t) =\int^1_0 f(\lambda,\theta,t)d\theta$.

In this notation, we have
$$h_1(\lambda,\theta,t)\in F^1\left(\frac{2n + 1}{n + 2}\right),\quad\quad h_2(\lambda,\theta,t) \in F^0\left(\frac{n + 1}{n + 2}\right).$$

In the following proposition we introduce a canonical transformation which transforms the Hamiltonian $h$ defined by (\ref{H1}) into new ones, which is more closer to an integrable Hamiltonian.

For this purpose, we consider a sightly general Hamiltonian system with the Hamiltonian
\begin{equation}\label{add1}
H(\lambda,\theta,t) = h_0(\lambda,t) + h_1(\lambda,\theta,t) + h_2(\lambda,\theta,t),
\end{equation}
where
$$
h_0\in F^1\left(\frac{2n+2}{n+2}\right),\quad h_1\in F^1(a_1),\quad h_2\in F^0(a_2)
$$
with
$$
D_{\lambda} h_0 \ge c\cdot\lambda^{\frac{2n+2}{n+2}-1}>0,\quad a_1\le\frac{2n+1}{n+2},\quad a_2=\frac{n+1}{n+2}.
$$

\begin{proposition}\label{Prop}
 There exists a canonical transformation which transforms the Hamiltonian $H$ defined by (\ref{add1}) into
\begin{equation}\label{H2}
\widetilde{H}(\mu,\phi,t) = \tilde{h}_0(\mu,t) + \tilde{h}_1(\mu,\phi,t) + \tilde{h}_2(\mu,\phi,t),
\end{equation}
where
$$
\tilde{h}_0\in F^1\left(\frac{2n + 2}{n + 2}\right),\quad \tilde{h}_1\in F^1\left(2a_1-\frac{2n+2}{n+2}\right),\quad \tilde{h}_2\in F^0\left(a_2\right),
$$
$$
\tilde{h}_0(\mu,t)= h_0(\mu,t)+[h_1](\mu,t).
$$
Moreover, the function $\tilde h_0(\mu,t)$ satisfies
$$
D_\mu\tilde h_0\ge c\cdot\mu^{\frac{2n+2}{n+2}-1}>0.
$$
\end{proposition}

\Proof We define the canonical transformation implicitly by
$$
\lambda= \mu + D_\theta S(\mu,\theta,t) =:\mu+ \nu, \quad\quad \phi = \theta + D_\mu S(\mu,\theta,t),
$$
where the generating function $S$ will be determined later. Under this transformation,
the new Hamiltonian is
$$
\begin{array}{lll}
\widetilde{H}(\mu,\theta,t) &=& h_0(\mu+ \nu, t) + h_1(\mu+ \nu, \theta, t) + h_2(\mu+ \nu, \theta, t)+ D_t S(\mu,\theta,t)\\[0.2cm]
&=& h_0(\mu,t) + D_\mu h_0(\mu,t)\nu + \int^1_0 (1-s)D^2_\mu h_0(\mu + s\nu,t)\nu^2ds\\[0.2cm]
&+& h_1(\mu, \theta, t)+\int^1_0 D_\mu h_1(\mu + s\nu,\theta, t)\nu ds\\[0.2cm]
&+& h_2(\mu+ \nu, \theta, t)+ D_t S(\mu,\theta,t).
\end{array}
$$

Let
$$
\nu(\mu,\theta,t)= \frac{h_1-[h_1]}{D_\mu h_0}\in F^1\left(a_1-\frac{n}{n+2}\right),
$$
and
$$
S(\mu,\theta,t)=\int_0^\theta\,\nu(\mu,s,t)ds.
$$
Then we obtain
$$
\begin{array}{lll}
\widetilde{H}(\mu,\theta,t) &=& h_0(\mu,t) + [h_1](\mu,t)\\[0.2cm]
&+& \int^1_0 (1-s)D^2_\mu h_0(\mu + s\nu,t)\nu^2ds+\int^1_0 D_\mu h_1(\mu + s\nu,\theta, t)\nu ds\\[0.2cm]
&+& h_2(\mu+ \nu, \theta, t)+ D_t S(\mu,\theta,t)\\[0.2cm]
&:=&\bar{h}_0(\mu,t)+\bar{h}_1(\mu,\theta, t)+\bar{h}_2(\mu, \theta, t),
\end{array}
$$
where
$$
\begin{array}{lll}
\bar{h}_0(\mu,t) &=& h_0(\mu,t) + [h_1](\mu,t),\\[0.2cm]
\bar{h}_1(\mu,\theta,t) &=& h_0(\mu+ \nu, t)-h_0(\mu,t) - D_\mu h_0(\mu,t)\nu+h_1(\mu+ \nu, \theta, t)-h_1(\mu, \theta, t),\\[0.2cm]
\bar{h}_2(\mu,\theta,t) &=&  h_2(\mu+ \nu, \theta, t)+ D_t S(\mu,\theta,t).
\end{array}
$$

It is easy to see that
$$
\bar{h}_0\in F^1\left(\frac{2n + 2}{n + 2}\right).
$$
 Moreover,
$$
\begin{array}{lll}
\bar{h}_1(\mu,\theta,t) &=& \dint_0^1(1-s)D_\mu^2h_0(\mu+s\nu,t)\nu^2ds + \dint_0^1D_\mu h_1(\mu+s\nu,\theta,t)\nu ds\\[0.42cm]
&\in& F^1\left(2a_1-\frac{2n+2}{n+2}\right)
\end{array}
$$
and
$$
\bar{h}_2(\mu,\theta,t)\in F^0(a_2).
$$

From the equation $ \phi= \theta+ D_\mu S(\mu,\theta,t)$, we can solve that $\theta = \phi+\Theta(\mu,\phi,t)$ with $\Theta \in F^1\left(a_1-\frac{2n+2}{n + 2}\right)$, and the new Hamiltonian is
$$
\begin{array}{lll}
\widetilde{H}(\mu,\phi,t) &=& \widetilde{H}(\mu,\phi+\Theta(\mu,\phi,t), t)\\[0.2cm]
&=& \bar{h}_0(\mu,t)+ \bar{h}_1(\mu,\phi+\Theta(\mu,\phi,t), t)+\bar{h}_2(\mu, \phi+\Theta(\mu,\phi,t), t)\\[0.2cm]
&:=&\tilde{h}_0(\mu,t) + \tilde{h}_1(\mu,\phi,t) + \tilde{h}_2(\mu,\phi,t),
\end{array}
$$
where
$$
\tilde{h}_0\in F^1\left(\frac{2n + 2}{n + 2}\right),\quad \tilde{h}_1\in F^1\left(2a_1-\frac{2n+2}{n+2}\right),\quad \tilde{h}_2\in F^0\left(a_2\right).
$$
Thus we have finished the proof of this proposition.\qed

Applying this proposition many times to the original system (\ref{H1}), the transformed Hamiltonian is of the form
\begin{equation}\label{add2}
H(\lambda,\theta,t) = h_0(\lambda,t) + h_1(\lambda,\theta,t) + h_2(\lambda,\theta,t),
\end{equation}
which has the property
$$
h_0\in F^1\left(\frac{2n + 2}{n + 2}\right),\quad h_1\in F^1\left(\frac{1-4n}{n + 2}\right),\quad h_2\in F^0\left(\frac{n + 1}{n + 2}\right),
$$
here we use the original notation.

Because $h_2\in F^0\left(\frac{n + 1}{n + 2}\right)$, let $\bar{h}_2(\lambda,\theta,t) = \lambda^{-\frac{n+1}{n+2}}h_2(\lambda,\theta,t)$, then there is a constant $M>0$ such that, for $\lambda\geq \lambda_0 ,(\theta,t)\in \mathbb{S}\times\mathbb{S}$
$$
|\lambda^jD^j_\lambda D^l_\theta\bar{h}_2(\lambda,\theta,t)|\le M,
$$
where $j+l$ is bounded by a sufficiently large integer $K>0$, which guarantees the smoothness assumption $C^4$ in the large twist theorem for  the Poincar\'{e} map (\ref{P}), hence we can choose the constant $M$ independently of $j,l$ for $j+l\leq K$.

For $j+l\leq K$, we expand the function $\bar{h}_{2}^{jl}(\lambda,\theta,t):=\lambda^jD^j_\lambda D^l_\theta\bar{h}_2(\lambda,\theta,t)$ into the Fourier series with respect to $t$
$$
\frac{\bar{h}_{2}^{jl0}}{2}+\sum_{n=1}^{+\infty}\left(\bar{h}_{2}^{jlnc}\cos 2\pi nt+\bar{h}_{2}^{jlns}\sin 2\pi nt\right)
$$
with the Fourier coefficients $\bar{h}_{2}^{jlnc}, \bar{h}_{2}^{jlns}$,  and its part sum is
$$
S_n(\bar{h}_{2}^{jl})=\frac{\bar{h}_{2}^{jl0}}{2}+\sum_{k=1}^{n}\left(\bar{h}_{2}^{jlkc}\cos 2\pi kt+\bar{h}_{2}^{jlks}\sin 2\pi kt\right).
$$

Consider the Fej\'{e}r sum of $\bar{h}_{2}^{jl}$
$$
F_n(\bar{h}_{2}^{jl})=\frac{\dsum_{k=0}^n S_k(\bar{h}_{2}^{jl})}{n+1}.
$$
By Fej\'{e}r Theorem,  the Fej\'{e}r sum $F_n(\bar{h}_{2}^{jl})$ of $\bar{h}_{2}^{jl}$ converges to $\bar{h}_{2}^{jl}$ uniformly with respect to $t$.
For any fixed $\epsilon > 0$, there exists a sufficiently large $N_{jl}\in \mathbb{N}$ such that for $n\geq N$ and $\lambda\geq \lambda_0 ,(\theta,t)\in \mathbb{S}\times\mathbb{S}$,
$$
|\bar{h}_{2}^{jl}(\lambda,\theta,t)-F_n(\bar{h}_{2}^{jl})(\lambda,\theta,t)|<\epsilon.
$$
Since the integers $j,l$ satisfying $j+l\leq K$ are finite, the positive integer $N_{jl}$ can be chosen independently of $j,l$. Let
$$
\bar{h}_{21}(\lambda,\theta,t)=F_N(\bar{h}_{2}^{00})(\lambda,\theta,t), \ \ \ \ \bar{h}_{22}(\lambda,\theta,t)=\bar{h}_{2}^{00}(\lambda,\theta,t)-F_N(\bar{h}_{2}^{00})(\lambda,\theta,t).
$$
Then
$$
\bar{h}_2(\lambda,\theta,t) = \bar{h}_{21}(\lambda,\theta,t) + \bar{h}_{22}(\lambda,\theta,t),
$$
and the functions $\bar{h}_{21}$ and $\bar{h}_{22}$ satisfy
 for $\lambda\geq \lambda_0 ,(\theta,t)\in \mathbb{S}\times\mathbb{S}$, and $j+1\leq K$
$$
|\lambda^jD^j_\lambda D^l_\theta\bar{h}_{21}(\lambda,\theta,t)|\le M+1,
$$
$$
|\lambda^jD^j_\lambda D^l_\theta\bar{h}_{22}(\lambda,\theta,t)|<\epsilon.
$$

From the above discussions, we may rewrite the Hamiltonian (\ref{add2}) in the form
\begin{equation}\label{add3}
H(\lambda,\theta,t) = h_0(\lambda,t) + h_1(\lambda,\theta,t) + h_{21}(\lambda,\theta,t) + h_{22}(\lambda,\theta,t),
\end{equation}
where
$$
h_0\in F^1\left(\frac{2n+2}{n+2}\right),\quad h_1\in F^1\left(\frac{1-4n}{n+2}\right),
$$
$$
h_{21}\in F^\infty\left(\frac{n+1}{n+2}\right), \quad h_{22}\in F^0\left(\frac{n+1}{n+2}\right),
$$
and
$$
\sup_{\lambda\ge\lambda_0,(\theta,t)\in\mathbb{S}\times\mathbb{S}}\left|\lambda^{j-\frac{n+1}{n+2}}D^j_\lambda D^l_\theta h_{22}(\lambda,\theta,t)\right|<\epsilon.
$$

Applying Proposition \ref{Prop} to the Hamiltonian (\ref{add3}), we may assume that the transformed hamiltonian of (\ref{add3}) is of the form
\begin{equation}\label{add4}
H(\lambda,\theta,t) = h_0(\lambda,t) + h_1(\lambda,\theta,t) + h_2(\lambda,\theta,t)
\end{equation}
where the functions $h_0, h_1$ and $h_2$ satisfy
$$
h_0(\lambda,t) = d\cdot\lambda^{\frac{2n+2}{n+2}} + \bar h_0(\lambda, t),
$$
and
$$
\bar h_0\in F^1\left(\frac{2n+1}{n+2}\right),\quad h_1\in F^1\left(\frac{1-4n}{n+2}\right),\quad h_2\in F^0\left(\frac{n+1}{n+2}\right).
$$
Moreover,
$$
\sup_{\lambda\ge\lambda_0,(\theta,t)\in\mathbb{S}\times\mathbb{S}}\left|\lambda^{j-\frac{n+1}{n+2}}D^j_\lambda D^l_\theta h_{2}(\lambda,\theta,t)\right|<\epsilon.
$$

The corresponding Hamiltonian system is
\begin{equation*}
\left\{\begin{array}{lll}
\dot{\lambda}&=&-D_{\theta}h_1(\lambda,\theta,t)-D_{\theta}h_{2}(\lambda,\theta,t),\\[0.2cm]
\dot{\theta}&=&D_{\lambda}h_0(\lambda,t)+D_{\lambda}h_1(\lambda,\theta,t)+D_{\lambda}h_{2}(\lambda,\theta,t).
 \end{array}\right.
\end{equation*}
Define
$$
\rho=\lambda^{\frac{1}{n+2}}, \quad\quad \theta=\theta,
$$
then
$$
\begin{array}{lll}
\frac{d\rho}{dt}&=&\frac{1}{n+2}\lambda^{-\frac{n+1}{n+2}}\frac{d\lambda}{dt}\\[0.2cm]
&=&-\frac{1}{n+2}\lambda^{-\frac{n+1}{n+2}}D_{\theta}h_1(\lambda,\theta,t)-\frac{1}{n+2}\lambda^{-\frac{n+1}{n+2}}D_{\theta}h_{2}(\lambda,\theta, t)\\[0.2cm]
&:=&f_{1}(\rho,\theta,t)+f_{2}(\rho,\theta, t),\\[0.2cm]
\frac{d\theta}{dt}&=&D_{\lambda}h_0(\lambda,t)+D_{\lambda}h_1(\lambda,\theta,t)+D_{\lambda}h_{2}(\lambda,\theta, t)\\[0.2cm]
&:=&f_{3}(\rho,t)+f_{4}(\rho,\theta,t)+f_{5}(\rho,\theta, t),
\end{array}
$$
that is,
\begin{equation*}
\left\{\begin{array}{lll}
\dot{\rho}&=&f_{1}(\rho,\theta,t)+f_{2}(\rho,\theta, t),\\[0.2cm]
\dot{\theta}&=&f_{3}(\rho,t)+f_{4}(\rho,\theta,t)+f_{5}(\rho,\theta, t),
 \end{array}\right.
\end{equation*}
where
$$
f_1\in F^1(-5n), f_2\in F^0(0),f_3\in F^1(n), f_4\in F^1(-5n-1), f_5\in F^0(-1).
$$
Moreover,
$$
f_3 = d\cdot\rho^n + O(\rho^{n-1}),\quad |\rho^jD^j_\rho D^l_\theta f_2|<\epsilon,\ \  |\rho^{j+1}D^j_\rho D^l_\theta f_5| <\epsilon.
$$

Now we introduce a parameter $\gamma$ by
$$
\rho=\gamma^{-1}I,\quad\quad I\in[1,2].
$$
Thus $\rho\to+\infty \Leftrightarrow \gamma\to 0^+$,  and the corresponding system has the form
\begin{equation}\label{P0}
\begin{array}{lll}
\dot{I}&=&\gamma f_{1}(\gamma^{-1}I,\theta,t)+\gamma f_{2}(\gamma^{-1}I,\theta, t),\\[0.2cm]
       &:=&\gamma^{5n+1} \bar{f}_{1}(I,\theta,t;\gamma)+\gamma \bar{f}_{2}(I,\theta,t;\gamma),\\[0.2cm]
\dot{\theta}&=&f_{3}(\gamma^{-1}I,t)+f_{4}(\gamma^{-1}I,\theta,t)+f_{5}(\gamma^{-1}I,\theta, t)\\[0.2cm]
&:=&\gamma^{-n}\bar{f}_{3}(I,t;\gamma)+\gamma^{5n+1}\bar{f}_{4}(I,\theta,t;\gamma)+\gamma \bar{f}_{5}(I,\theta, t;\gamma),
 \end{array}
\end{equation}
where $\bar{f}_i$ $(i=1,2,3,4,5)$ satisfy
\begin{equation}\label{A0}
|D_I^k\bar{f}_3(I,t;\gamma)|\leq C,
\end{equation}
\begin{equation}\label{A1}
|D_I^k D_{\theta}^j \bar{f}_1(I,\theta,t;\gamma)|\leq C,\ \  |D_I^k D_{\theta}^j \bar{f}_4(I,\theta,t;\gamma)|\leq C
\end{equation}
\begin{equation}\label{A2}
|D_I^k D_{\theta}^j \bar{f}_2(I,\theta,t;\gamma)| < \epsilon,\ \  |D_I^k D_{\theta}^j \bar{f}_5(I,\theta,t;\gamma)| < \epsilon.
\end{equation}
for some constant $C>0$ and $k+j\leq 4$.

Let $(I(t),\theta(t))$ be the solution of (\ref{P0}) with the initial condition $(I(t),\theta(t))=(I_0, \theta_0)$, then by the standard method (for example, see \cite{Diecherhoff}, \cite{Wang96}), we know that the Poincar\'{e} map has the form
\begin{equation}\label{P}
\left\{\begin{array}{lll}
I_1&=&I_0+  \Xi_1(\theta_0, I_0; \gamma),\\[0.2cm]
\theta_1&=&\theta_0+\gamma^{-n}\psi(I_0;\gamma)+ \Xi_2(\theta_0, I_0; \gamma),
 \end{array}\right.
\end{equation}
where
$\psi'(I_0;\gamma)\geq \frac{d}{2}$ for $I_0\in [1,2]$ with $\gamma$ sufficiently small, and
\begin{equation}\label{Small1}
|D_{I_0}^k D_{\theta_0}^l\Xi_i(\theta_0, I_0; \gamma)|\leq O(\gamma^{n+1})+\gamma^{-4n}\epsilon \,O(\gamma)
\end{equation}
as $\gamma\to 0^+$ for $i=1,2$ and $k+l\leq 4$.

Indeed, let
$$
\begin{array}{lll}
B_1(I_0,\theta_0,t;\gamma)&=&\int_0^t\bar{f}_{3}(I(s),s;\gamma)ds,\\[0.2cm]
A(I_0,\theta_0,t;\gamma)&=&\gamma^{5n+1} \int_0^t\bar{f}_{1}(I(s),\theta(s),s;\gamma)ds+\gamma \int_0^t\bar{f}_{2}(I(s),\theta(s),s;\gamma)ds,\\[0.2cm]
B_2(I_0,\theta_0,t;\gamma)&=&\gamma^{5n+1} \int_0^t\bar{f}_{4}(I(s),\theta(s),s;\gamma)ds+\gamma \int_0^t\bar{f}_{5}(I(s),\theta(s),s;\gamma)ds,
\end{array}
$$
then
$$
\begin{array}{lll}
I(t)&=&I_0+A(I_0,\theta_0,t;\gamma),\\[0.2cm]
\theta(t)&=&\theta_0+\gamma^{-n} B_1(I_0,\theta_0,t;\gamma)+B_2(I_0,\theta_0,t;\gamma)
\end{array}
$$
and functions $A,B_1,B_2$ satisfy
\begin{equation}\label{ABC}
\begin{array}{lll}
B_1&=&\int_0^t\bar{f}_{3}(I_0+A,s;\gamma)ds\\[0.2cm]
A&=&\gamma^{5n+1} \int_0^t\bar{f}_{1}(I_0+A,\theta_0+\gamma^{-n}B_1+B_2,s;\gamma)ds\\[0.2cm]
&+&\gamma \int_0^t\bar{f}_{2}(I_0+A,\theta_0+\gamma^{-n}B_1+B_2,s;\gamma)ds,\\[0.2cm]
B_2&=&\gamma^{5n+1} \int_0^t\bar{f}_{4}(I_0+A,\theta_0+\gamma^{-n}B_1+B_2,s;\gamma)ds\\[0.2cm]&+&\gamma \int_0^t\bar{f}_{5}(I_0+A,\theta_0+\gamma^{-n}B_1+B_2,s;\gamma)ds.
\end{array}
\end{equation}
Therefore
$$
\begin{array}{lll}
\psi(I_0;\gamma)&=&\int_0^1\,\bar{f}_{3}(I(t),t;\gamma)dt\\[0.2cm]
&=&\int_0^1\,[d \cdot I(t)^n+ I(t)^{n-1}O(\gamma)]dt\\[0.2cm]
&=&d \cdot (I_0+O(\gamma^{5n+1})+\epsilon O(\gamma))^n+(I_0+O(\gamma^{5n+1})+\epsilon O(\gamma))^{n-1}O(\gamma )\\[0.2cm]
&=&d \cdot I_0^n+O(\gamma)(1+\epsilon)I_0^{n-1}+\cdots,\\[0.2cm]
\end{array}
$$
which implies that
$$
\psi'(I_0;\gamma)\geq \frac{d}{2}
$$
for $I_0\in [1,2]$ with $\gamma$ sufficiently small, and
$$
\begin{array}{lll}
 \Xi_1(\theta_0, I_0; \gamma)&=&\gamma^{5n+1} \int_0^1\bar{f}_{1}(I_0+A,\theta_0+\gamma^{-n}B_1+B_2,s;\gamma)ds\\[0.2cm]&+&\gamma \int_0^1\bar{f}_{2}(I_0+A,\theta_0+\gamma^{-n}B_1+B_2,s;\gamma)ds,\\[0.2cm]
 \Xi_2(\theta_0, I_0; \gamma)&=&\gamma^{5n+1} \int_0^1\bar{f}_{4}(I_0+A,\theta_0+\gamma^{-n}B_1+B_2,s;\gamma)ds\\[0.2cm]&+&\gamma \int_0^1\bar{f}_{5}(I_0+A,\theta_0+\gamma^{-n}B_1+B_2,s;\gamma)ds.
\end{array}
$$
Firstly, according to (\ref{A1}) and (\ref{A2}), it is easy to see that
$$|\Xi_i(\theta_0, I_0; \gamma)|\leq O(\gamma^{5n+1})+\epsilon \,O(\gamma)$$
as $\gamma\to 0^+$ for $i=1,2$. As for the estimates on the derivatives, from the equations (\ref{ABC}), we find that
\begin{equation*}
\begin{array}{lll}
D_{I_0}B_1&=&\int_0^tD_{I}\bar{f}_{3}\cdot(1+D_{I_0}A)ds\\[0.2cm]
D_{I_0}A&=&\gamma^{5n+1} \int_0^t\left[D_{I}\bar{f}_{1}\cdot(1+D_{I_0}A)+D_{\theta}\bar{f}_{1}\cdot(\gamma^{-n}D_{I_0}B_1+D_{I_0}B_2)\right]ds\\[0.2cm]
&+&\gamma \int_0^t\left[D_{I}\bar{f}_{2}\cdot(1+D_{I_0}A)+D_{\theta}\bar{f}_{2}\cdot(\gamma^{-n}D_{I_0}B_1+D_{I_0}B_2)\right]ds,\\[0.2cm]
D_{I_0}B_2&=&\gamma^{5n+1} \int_0^t\left[D_{I}\bar{f}_{4}\cdot(1+D_{I_0}A)+D_{\theta}\bar{f}_{4}\cdot(\gamma^{-n}D_{I_0}B_1+D_{I_0}B_2)\right]ds\\[0.2cm]
&+&\gamma \int_0^t\left[D_{I}\bar{f}_{5}\cdot(1+D_{I_0}A)+D_{\theta}\bar{f}_{5}\cdot(\gamma^{-n}D_{I_0}B_1+D_{I_0}B_2)\right]ds,
\end{array}
\end{equation*}
which together with (\ref{A0}), (\ref{A1}),(\ref{A2}) implies that
$$|D_{I_0}\Xi_i(\theta_0, I_0; \gamma)|\leq O(\gamma^{4n+1})+\gamma^{-n}\epsilon \,O(\gamma)$$
for $i=1,2$. Differentiating $A,B,C$ with respect to $I_0$ $k$ times yields that
$$|D_{I_0}^k\Xi_i(\theta_0, I_0; \gamma)|\leq O(\gamma^{5n+1-kn})+\gamma^{-kn}\epsilon \,O(\gamma)$$
for $i=1,2$ and $0\leq k\leq 4$. Also by Differentiating $A,B,C$ with respect to $\theta_0$, we know that
\begin{equation*}
\begin{array}{lll}
D_{\theta_0}B_1&=&\int_0^tD_{I}\bar{f}_{3}\cdot D_{\theta_0}Ads\\[0.2cm]
D_{\theta_0}A&=&\gamma^{5n+1} \int_0^t\left[D_{I}\bar{f}_{1}\cdot D_{\theta_0}A+D_{\theta}\bar{f}_{1}\cdot(1+\gamma^{-n}D_{\theta_0}B_1+D_{\theta_0}B_2)\right]ds\\[0.2cm]
&&+\gamma \int_0^t\left[D_{I}\bar{f}_{2}\cdot D_{\theta_0}A+D_{\theta}\bar{f}_{2}\cdot(1+\gamma^{-n}D_{\theta_0}B_1+D_{\theta_0}B_2)\right]ds,\\[0.2cm]
D_{\theta_0}B_2&=&\gamma^{5n+1} \int_0^t\left[D_{I}\bar{f}_{4}\cdot D_{\theta_0}A+D_{\theta}\bar{f}_{4}\cdot(1+\gamma^{-n}D_{\theta_0}B_1+D_{\theta_0}B_2)\right]ds\\[0.2cm]
&&+\gamma \int_0^t\left[D_{I}\bar{f}_{5}\cdot D_{\theta_0}A+D_{\theta}\bar{f}_{5}\cdot(1+\gamma^{-n}D_{\theta_0}B_1+D_{\theta_0}B_2)\right]ds,
\end{array}
\end{equation*}
which together with (\ref{A0}), (\ref{A1}),(\ref{A2}) implies that
$$|D_{\theta_0}\Xi_i(\theta_0, I_0; \gamma)|\leq O(\gamma^{5n+1})+\epsilon \,O(\gamma), \ \ i=1,2.$$
The estimates on the higher derivatives $D_{I_0}^k D_{\theta_0}^l\Xi_i$ can be proven similarly, and thus we have finished the proof of (\ref{Small1}).

Now we can choose $\epsilon=\gamma^{5n}$ with $\gamma$ small enough such that the functions $\Xi_1,\Xi_2$ satisfy the smallness condition in Theorem \ref{A},  then this theorem can be applied, and one can obtain the boundedness of solutions and the existence of quasi-periodic solutions stated in the main result as usual. This finishes the proof of Theorem \ref{B}.

\bigskip

\end{document}